\documentclass[11pt]{amsart}

\usepackage{amsfonts,amssymb,stmaryrd,amscd,amsmath,latexsym,amsbsy}

\usepackage{amssymb}
\usepackage{amsfonts}
\usepackage{latexsym}

\newtheorem{theorem}{Theorem}[section]

\newtheorem{corollary}[theorem]{Corollary}
\theoremstyle{definition}

\newtheorem{remark}[theorem]{Remark}

\newcommand{\g}{\mathfrak{g}}
\newcommand{\h}{\mathfrak{h}}

\newcommand{\ben}{\begin{enumerate}}
\newcommand{\een}{\end{enumerate}}

\theoremstyle{plain}

\newtheorem*{sol}{Solution}

\theoremstyle{definition}

\theoremstyle{remark}

\newcommand{\solu}[1]{\begin{sol}{\bf (\ref{#1})}}

\begin{document}

\title{Reducibility of the polynomial representation 
of the degenerate double affine Hecke algebra}

\author{Pavel Etingof}
\address{Department of Mathematics, Massachusetts Institute of Technology,
Cambridge, MA 02139, USA}
\email{etingof@math.mit.edu}

\maketitle

\section{Introduction}

In this note we determine the values of parameters $c$ for which the
polynomial representation of the degenerate double affine Hecke
algebra (DAHA), i.e. the trigonometric Cherednik algebra, is
reducible. Namely, we show that $c$ is a reducibility point 
for the polynomial representation of the trigonometric 
Cherednik algebra for a root system $R$ if and only if it is a 
reducibility point for the {\it rational}
Cherednik algebra for the Weyl group of some root subsystem $R'\subset R$
of the same rank; such subsystems for any $R$ are given by the
well known Borel-de Siebenthal algorithm.   

This generalizes to the trigonometric case the 
result of \cite{DJO}, where the reducibility points are found for
the rational Cherednik algebra. Together with the result of
\cite{DJO}, our result gives an explicit list of reducibility 
points in the trigonometric case. 

We emphasize that our result is contained in the recent work 
of I. Cherednik \cite{Ch2}, where reducibility points 
are determined for nondegenerate DAHA. Namely, the techniques 
of \cite{Ch2}, based on intertwiners, work equally well in 
the degenerate case. In fact, outside of roots of unity, 
the questions of reducibility of the polynomial representation 
for the degenerate and nondegenerate DAHA are equivalent, 
and thus our result is equivalent to that of 
\cite{Ch2}. However, our proof is quite different from that in
\cite{Ch2}; it is based on the geometric approach 
to Cherednik algebras developed in \cite{E2},
and thus clarifies the results of \cite{Ch2} from a geometric
point of view. In particular, we explain that our result 
and its proof can be generalized 
to the much more general setting of Cherednik algebras for  
any smooth variety with a group action.   
 
We note that in the non-simply laced case, it is not true 
that the reducibility points for $R$ are the same in the
trigonometric and rational settings. In the trigonometric
setting, one gets additional reducibility points, which 
arise for type $B_n$, $n\ge 3$, $F_4$, and
$G_2$, but not for $C_n$. This phenomenon was discovered by
Cherednik (in the $B_n$ case, see \cite{Ch3}, Section 5);
in \cite{Ch2}, he gives a complete list of additional
reducibility points. At first sight, this list looks somewhat 
mysterious; our work demystifies it, by interpreting it in 
terms of the Borel - de Siebenthal classification of 
equal rank embeddings of root systems. 

The result of this note is a manifestation of the general
principle that the representation theory of the trigonometric
Cherednik algebra (degenerate DAHA) for a root system $R$ reduces
to the representation theory of the {\it rational} Cherednik 
algebra for Weyl groups of root subsystems $R'\subset R$. 
This principle is the ``double'' analog of a similar
principle in the representation theory of affine Hecke algebras, 
which goes back to the work of Lusztig \cite{L}, in which it is
shown that irreducible representations of the affine Hecke algebra 
of a root system $R$ may be described in terms of irreducible 
representations of 
the degenerate affine Hecke algebras for Weyl groups of root subsystems $R'\subset
R$. We illustrate this principle at the end of the note by applying it to finite
dimensional representations of trigonometric Cherednik algebras.  

{\bf Acknowledgements.} The author is very grateful to
I. Cherednik for many useful discussions, and for sharing the
results of his work \cite{Ch2} before its publication. 
The author also thanks G. Lusztig, M. Varagnolo, and E. Vasserot for useful discussions. 
The work of the author was  partially supported by the NSF grant
DMS-0504847.

\section{Preliminaries}

\subsection{Preliminaries on root systems}\label{prel}

Let $W$ be a irreducible Weyl group, $\h$ its (complex)
reflection representation, and $L\subset \h$ a $\Bbb Z$-lattice invariant
under $W$. 

For each reflection $s\in W$, let $L_s$ be the intersection of
$L$ with the $-1$-eigenspace of $s$ in $\h$, and let
$\alpha_s^\vee$ be a generator of $L_s$. Let $\alpha_s$ be the
element in $\h^*$ such that $s\alpha_s=-\alpha_s$, and
$(\alpha_s,\alpha_s^\vee)=2$. Then we have
$$
s(x)=x-(x,\alpha_s)\alpha_s^\vee,\ x\in \h.
$$

Let $R\subset \h^*$ be the collection of vectors $\pm \alpha_s$,
and $R^\vee\subset \h$ the collection of vectors
$\pm\alpha_s^\vee$. It is well known that $R,R^\vee$ are mutually
dual reduced root systems.  Moreover, we have $Q^\vee\subset
L\subset P^\vee$, where $P^\vee$ is the coweight lattice, and
$Q^\vee$ the coroot lattice.

Consider the simple complex Lie group $G$ with root system $R$,
whose center is $P^\vee/L$. The maximal torus of $G$ can be
identified with $H=\h/L$ via the exponential map. 

For $g\in H$, let $C_\g(g)$ be the centralizer of $g$ in
$\g:={\rm Lie}(G)$.  Then $C_\g(g)$ is a reductive subalgebra of
$\g$ containing $\h$, and its Weyl group is the stabilizer $W_g$
of $g$ in $W$.

Let $\Sigma\subset H$ be the set of elements whose centralizer
$C_\g(g)$ is semisimple (of the same rank as $\g$). 
$\Sigma$ can also be defined as the set of point strata for the
stratification of $H$ with respect to stabilizers. It is well
known that the set $\Sigma$ is finite, and the Dynkin diagram of
$C_\g(g)$ is obtained from the extended Dynkin diagram of $\g$ by
deleting one vertex (the Borel-de Siebenthal
algorithm). Moreover, any Dynkin diagram obtained in this way
corresponds to $C_\g(g)$ for some $g$.

\subsection{The rational Cherednik algebra}

Let $W$ be a Coxeter group with reflection representation
$\h$. For any conjugacy invariant function $c$ on the set of
reflections in $W$, one can define the rational Cherednik algebra
$H_{1,c}(W,\h)$ (see e.g. \cite{E1}, Section 7); we will denote it shortly by
$H_c(W,\h)$. This algebra has a polynomial representation $\Bbb
C[\h]$, which is defined using Dunkl operators. A function $c$ is
said to be {\it singular} if the polynomial representation is
reducible. Let ${\rm Sing}(W,\h)$ be the set of singular $c$. This
set is determined explicitly in \cite{DJO}.

\subsection{The degenerate DAHA}

Let $W,L,H$ be as in subsection \ref{prel}. A reflection
hypertorus in $H$ is a connected component $T$ of the fixed set
$H^s$ for a reflection $s\in W$. Let $c$ be a conjugation
invariant function on the set of reflection hypertori. 

The degenerate DAHA attached to $W,H$ was introduced by Cherednik,
see e.g. \cite{Ch1}. This algebra is generated by polynomial functions on
$H$, the group $W$, and trigonometric Dunkl operators.
Using the geometric approach of \cite{E2}, which attaches a 
Cherednik algebra to any smooth affine algebraic variety with a 
finite group action, the degenerate DAHA can be defined as 
the Cherednik algebra $H_{1,c}(W,H)$ attached to the variety $H$
with the action of the finite group $W$. We will shortly denote
this algebra by $H_c(W,H)$.  

Note that this setting includes the case of non-reduced root systems. 
Namely, in the case of a non-reduced root system the function $c$   
may take nonzero values on reflection hypertori which 
don't go through $1\in H$. 

\section{The results} 

\subsection{The main results} 

The degenerate DAHA has a polynomial representation $M=\Bbb C[H]$
on the space of regular functions on $H$. We would like to
determine for which $c$ this representation is reducible.

Let $g\in \Sigma$. Denote by $c_g$ the restriction of the
function $c$ to reflections in $W_g$; that is, for $s\in W_g$,
$c_g(s)$ is the value of $c$ on the (unique) hypertorus 
$T_{g,s}$ passing through $g$ and fixed by $s$. Denote by ${\rm
Sing}_g(W,L)$ the set of $c$ such that $c_g\in {\rm Sing}(W_g,\h)$.

\begin{remark}
If $c(T)=0$ unless $T$ contains $1\in H$ (``reduced
case'') then $c$ can be regarded as a function of
reflections in $W$, and $c_g$ is the usual restriction of $c$ to
reflections in $W_g$. 
\end{remark} 

Our main result is the following. 

\begin{theorem}\label{irr}
The polynomial representation $M$ of $H_c(W,H)$ is reducible if
and only if $c\in \cup_{g\in \Sigma}{\rm Sing}_g(W,L)$. 
\end{theorem} 

The proof of this theorem is given in the next subsection. 

\begin{corollary}\label{irr1} 
If $c$ is a constant function (in particular, if $R$ is simply
laced), then the polynomial representation $M$ of $H_c(W,H)$ 
is reducible if and only if so is the polynomial representation
of the rational Cherednik algebra $H_c(W,\h)$, i.e. iff 
$c=j/d_i$, where $d_i$ is a degree of $W$, and $j$ is a positive
integer not divisible by $d_i$.  
\end{corollary}

\begin{proof}
The result follows from Theorem \ref{irr}, the result of
\cite{DJO}, and the well known fact\footnote{This fact is proved as follows. 
Let $P_W(t)$ be the Poincar\'e polynomial of $W$; 
so 
$$
P_W(t)=\prod_i \frac{1-t^{d_i(W)}}{1-t},
$$  
where $d_i(W)$ are the degrees of $W$. 
Then by Chevalley's theorem, $P_W(t)/P_{W'}(t)$ is a polynomial
(the Hilbert polynomial of the generators of 
the free module $\Bbb C[\h]^{W'}$ over $\Bbb C[\h]^W$).  
So, since the denominator vanishes at a root of unity
of degree $d_i(W')$, so does the numerator, which implies the
statement.} that for any subgroup
$W'\subset W$ generated by reflections, every degree of $W'$
divides some degree of $W$. 
\end{proof} 

However, if $c$ is not a constant function, 
the answer in the trigonometric case may differ from 
the rational case, as explained below.  

\subsection{Proof of Theorem \ref{irr}}

Assume first that the polynomial representation $M$ is
reducible. Then there exists a nonzero proper submodule $I\subset M$,
which is an ideal in $\Bbb C[H]$. This ideal defines a subvariety
$Z\subset H$, which is $W$-invariant; it is the support 
of the module $M/I$. It is easy to show using the results of
\cite{E2} (see e.g. \cite{BE}) that $Z$ is a union of strata 
of the stratification of $H$ with respect to stabilizers. 
In particular, since $Z$ is closed, it contains a stratum which
consists of one point $g$. Thus $g\in \Sigma$. 
Consider the formal completion $\widehat{M}_g$ of $M$ at $g$. 
As follows from \cite{E2} (see also \cite{BE}), 
this module can be viewed as a module over the formal completion
$\widehat{H}_{c_g}(W_g,\h)_0$ of the rational Cherednik algebra
of the group $W_g$ at $0$, and it has a nonzero proper  
submodule $\widehat{I}_g$. Thus, $\widehat{M}_g$ is reducible,
which implies (by taking nilpotent vectors under $\h^*$)
that the polynomial representation $\bar M$ over
$H_{c_g}(W_g,\h)$ is reducible, hence $c_g\in {\rm Sing}(W_g)$,
and $c\in {\rm Sing}_g(W,L)$. 

Conversely, assume that $c\in {\rm Sing}_g(W,L)$, and thus 
$c_g\in {\rm Sing}(W_g)$. Then the polynomial representation
$\bar M$ of $H_{c_g}(W_g,\h)$ is reducible. 
This implies that the completion $\widehat{M}_g=\widehat{\Bbb
C[H]}_g$ is a reducible 
module over $\widehat{H}_{c_g}(W_g,\h)_0$, i.e. it contains a
nonzero proper submodule (=ideal) $J$. Let $I\subset \Bbb C[H]$
be the intersection of $\Bbb C[H]$ with $J$. Clearly, 
$I\subset M$ is a proper submodule (it does not contain $1$). 
So it remains to show that it is nonzero. 
To do so, denote by $\Delta$ a regular function on $H$
which has simple zeros on all the reflection hypertori.  
Then clearly $\Delta^n\in J$ for large enough $n$, 
so $\Delta^n\in I$. Thus $I\ne 0$ and the theorem is proved. 

\subsection{Reducibility points in the non-simply laced case}

In this subsection we will consider the reduced non-simply laced
case, i.e. the case of root systems of type $B_n,C_n$, $F_4$, and
$G_2$. In this case, $c$ is determined by two numbers $k_1$ and
$k_2$, the values of $c$ on reflections for long and short roots,
respectively. 

The set ${\rm Sing}(W,\h)$ 
is determined for these cases in \cite{DJO}, as the union of the
following lines (where $l\ge 1$, $u=k_1+k_2$, and $i=1,2$). 

$B_n=C_n$: 
$$
2jk_1+2k_2=l,\ l\ne 0\text{ mod }2, j=0,...,n-1,
$$
and 
$$
jk_1=l,\ (l,j)=1,\ j=2,...,n. 
$$

$F_4$: 
$$
2k_i=l,\ 2k_i+2u=l,\ l\ne 0\text{ mod }2; 
3k_i=l,\ l\ne 0\text{ mod }3; 
$$
$$
2u=l, 4u=l,\ l\ne 0\text{ mod }2; 
$$
$$
6u=l,\ l=1,5,7,11\text{ mod }12.
$$

$G_2$: 
$$
2k_i=l,\ l\ne 0\text{ mod }2;\ 3u=l,\ l\ne 0\text{ mod }3.
$$

By using Theorem \ref{irr}, we determine that 
the polynomial representation in the trigonometric case is
reducible on these lines and also on the following additional
lines: 

$B_n$, $n\ge 3$:
$$
(2p-1)k_1=2q, n/2<q\le n-1, p\ge 1, (2p-1,q)=1.
$$

$F_4$:
$$
6k_1+2k_2=l, 4k_1=l,\ l\ne 0\text{ mod }2. 
$$

$G_2$: 
$$
3k_1=l,\ l\ne 0\text{ mod }3. 
$$

In the $C_n$ case, we get no additional lines. 

Note that exactly the same list of additional reducibility points
appears in \cite{Ch2}.

\begin{remark} As explained above, 
the additional lines appear from particular equal rank
embeddings of root systems. Namely, the additional lines for $B_n$ appear from 
the inclusion $D_n\subset B_n$. 
The two series of additional lines for $F_4$ appear from the embeddings
$B_4\subset F_4$ and $A_3\times A_1\subset F_4$, respectively. 
Finally, the additional lines for $G_2$ appear from the embedding
$A_2\subset G_2$.
\end{remark}

\subsection{Generalizations}

Theorem \ref{irr} can be generalized, with essentially the same proof, to the
setting of any smooth variety with a group action, as defined in
\cite{E2}. 

Namely, let $X$ be a smooth algebraic variety, and $G$
a finite group acting faithfully on $X$. Let $c$ be a
conjugation invariant function on the set of pairs $(g,Y)$,
where $g\in G$, and $Y$ is a connected component of $X^g$ which
has codimension 1 in $X$. Let $H_{1,c,0,X,G}$ be
the corresponding sheaf of Cherednik algebras defined in
\cite{E2}. We have the polynomial representation ${\mathcal O}_X$
of this sheaf. 

Let $\Sigma\in X$ be the set of points with maximal stabilizer,
i.e. points whose stabilizer is bigger than that of nearby
points. Then $\Sigma$ is a finite set. For $x\in X$, let $G_x$ be the
stabilizer of $x$ in $G$; it is a finite subgroup of $GL(T_xX)$. 
Let $c_x$ be the function of reflections in $G_x$ defined by 
$c_x(g)=c(g,Y)$, where $Y$ is the reflection hypersurface 
passing through $x$ and fixed by $g$ pointwise. 
Let ${\rm Sing}_x(G,X)$ be the set of $c$ such that 
$c_x\in {\rm Sing}(G_x,T_xX)$ (where ${\rm Sing}(G_x,T_xX)$
denotes the set of values of parameters $c$ for which 
the polynomial representation of the rational Cherednik
algebra $H_c(G_x,T_xX)$ is reducible). 

Then we have the following theorem, whose statement and proof are direct
generalizations of those of Theorem \ref{irr} (which is obtained
when $G$ is a Weyl group and $X$ a torus). 

\begin{theorem}\label{irr2}
The polynomial representation ${\mathcal O}_X$ of $H_{1,c,0,X,G}$
is reducible if and only if $c\in \cup_{x\in \Sigma}{\rm
Sing}_x(G,X)$.  
\end{theorem}

Note that this result generalizes in a straightforward way to the case
when $X$ is a complex analytic manifold, and $G$ a discrete
group of holomorphic transformations of $X$. 

\subsection{Finite dimensional representations of the
degenerate double affine Heclke algebra}\footnote{The contents of
this subsection 
arose from a discussion of the author with M. Varagnolo and E. Vasserot.}

Another application of the approach of this note is a description
of the category of finite dimensional representations of the 
degenerate DAHA in terms of categories of finite dimensional
representations of rational Cherednik algebras. 
Namely, let $FD(A)$ denote the category of finite dimensional
representations of an algebra (or sheaf of algebras) $A$. 
Then in the setting of the previous subsection we have the
following theorem (see also Proposition 2.22 of \cite{E2}). 

Let $\Sigma'$ be a set of representatives of $\Sigma/G$ in
$\Sigma$.

\begin{theorem}\label{fd}
One has 
$$
FD(H_{1,c,0,X,G})=\oplus_{x\in \Sigma'}FD(H_{c_x}(G_x,T_xX)).
$$
\end{theorem}

\begin{proof}
Suppose $V$ is a finite dimensional representation of 
$H_{1,c,0,X,G}$. Then the support of $V$ is a union of finitely
many points, and these points must be strata of the
stratification of $X$ with respect to stabilizers, so 
they belong to $\Sigma$. This implies that $V=\oplus_{\xi\in
\Sigma/G} V_\xi$, where $V_\xi$ is supported on the orbit $\xi$. 
Taking completion of the Cherednik algebra at $\xi$, we 
can regard the fiber $(V_\xi)_x$ for $x\in \xi$ as a module over the rational
Cherednik algebra $H_{c_x}(G_x,T_xX)$ (see \cite{E2,BE}). In this way, $V$ gives rise to an
object of $\oplus_{x\in \Sigma'}FD(H_{c_x}(G_x,T_xX))$.

This procedure can be reversed; this implies the theorem. 
\end{proof} 

\begin{corollary}\label{trig}
One has 
$$
FD(H_c(W,H))=\oplus_{g\in \Sigma/W}FD(H_{c_g}(W_g,\h)).
$$
\end{corollary}

\begin{remark}
Recall that a representation of $H_c(W,H)$ is said to be
spherical if it is a quotient of the polynomial representation. 
It is clear that the categorical equivalence of Corollary
\ref{trig} preserves sphericity of representations (in both
directions). This implies that the results of the paper
\cite{VV}, which classifies spherical finite-dimensional
representations of the rational Cherednik algebras, in fact
yields, through Corollary \ref{trig}, the classification 
of spherical finite dimensional representations of 
degenerate DAHA, and hence of nondgenerate DAHA outside of roots
of unity.  
\end{remark}

\end{document}